\documentclass[11pt]{article}
\usepackage{amsmath,amsfonts,amssymb,mathrsfs}
\usepackage[mathscr]{eucal}

\usepackage{mathtools}                                                      
\mathtoolsset{showonlyrefs=true}                                    

\usepackage{eufrak}                                             

\usepackage{amsmath}                                                        
\usepackage{graphicx}                                                       
\usepackage{subfigure}
\usepackage{enumitem}                                                    
\usepackage{hyperref}
\usepackage{amssymb}                                                        
\usepackage{cancel}                                                             
\usepackage[normalem]{ulem}                                                                 
\usepackage{pstricks}
\usepackage{rotating}
\usepackage{lscape}
\usepackage[paperwidth=8.5in,paperheight=11in,top=1.25in, bottom=1.25in, left=1.00in, right=1.00in]{geometry}
\usepackage{mathtools}                                                      
\mathtoolsset{showonlyrefs=true}                                    
\usepackage{fixltx2e,amsmath}                                           
\MakeRobust{\eqref}

\usepackage{amsthm}                                                             
\allowdisplaybreaks                                                             
\usepackage{chngcntr}

\usepackage{titlesec}
\usepackage[titletoc,toc,title]{appendix}

\usepackage{apptools, etoolbox}

\AtAppendix{\counterwithin{theorem}{section}}

\usepackage{color}

\newtheorem{theorem}{Theorem}[section]

\newtheorem{corollary}[theorem]{Corollary}
\newtheorem{lemma}[theorem]{Lemma}
\newtheorem{definition}[theorem]{Definition}

\newtheorem{remark}[theorem]{Remark}
\newtheorem{assumption}[theorem]{Assumption}
\newtheorem{notation}[theorem]{Notation}

\def \c {{\chi}}

\def \costc {c_{0}}
\def \costd {M}
\def \coste {c_{2}}

\def \cca {c_{3}}
\def \ccb {c_{4}}
\def \ccc {c_{5}}
\def \ccd {C}

\def \ccf {c_{8}}
\def \ccg {c_{9}}
\def \cch {c_{10}}
\def \cci {c_{11}}
\def \ccl {C}

\def \Kol {\mathcal{K}_{\costd,B}}
\def \costv {\|v\|_{L^{2}([t,T])}^{2}}

\def \b {{\beta}}
\def \d {{\delta}}
\def \l {{\lambda}}
\def \L {{\Lambda}}
\def \G {{\Gamma}}
\def \s {{\sigma}}

\def \R {{\mathbb {R}}}

\def \x {{\xi}}

\def \t {{\tau}}
\def \t {{\tau}}
\def \n {{\nu}}
\def \m {{\mu}}

\def \z {{\zeta}}

\def \LL {{\L^1}}

\def \Gg {{\G^1}}

\def \g {{\gamma}}
\def \O {{\Omega}}
\def \phi {{\varphi}}
\def \div {{\text{\rm div}}}

\def\p{\partial}
\def \P {{{P}}}

\def \rdd {{\mathbb {R}}^{d+1}}
\def \rd {{\mathbb {R}}^{d}}

\def \c {{\chi}}

\def \b {{\beta}}
\def \d {{\delta}}
\def \l {{\lambda}}
\def \L {{\Lambda}}
\def \G {{\Gamma}}
\def \s {{\sigma}}

\def \Dil {\mathcal{D}}

\def \R  {{\mathbb {R}}}
\def \x {{\xi}}
\def \g {{\gamma}}

\def \t {{\tau}}
\def \n {{\nu}}
\def \m {{\mu}}

\def \z {{\zeta}}
\def \p {{\partial}}

\def \O {{\Omega}}

\def \d {{\delta}}

\def \c {{\chi}}

\def \b {{\beta}}
\def \d {{\delta}}

\def \G {\Ga}

\def \Ga {{\Gamma}}

\def \s {{\sigma}}

\def \R {{\mathbb {R}}}

\def \x {{\xi}}

\def \t {{\tau}}
\def \t {{\tau}}
\def \n {{\nu}}
\def \v {{\nu}}
\def \m {{\mu}}

\def \z {{\zeta}}

\def \g {{\gamma}}
\def \O {{\Omega}}
\def \phi {{\varphi}}
\def \div {{\text{\rm div}}}

\def\l {\lambda}

\def \P {{P}}

\def \L {\mathscr{L}}

\def \à {{\`a }}
\def \è {{\`e }}
\def \ò {{\`o }}
\def \ù {{\`u }}

\def\l {\varrho}

\makeatletter
\def\section{\@startsection {section}{1}{\z@}{3.25ex plus 1ex minus
 .2ex}{1.5ex plus .2ex}{\large\bf}}
\def\subsection{\@startsection{subsection}{2}{\z@}{3.25ex plus 1ex minus
 .2ex}{1.5ex plus .2ex}{\normalsize\bf}}
\@addtoreset{equation}{section} 
\makeatother

\begin{document}

\title{Gaussian lower bounds for non-homogeneous Kolmogorov equations with measurable coefficients}

\author{Alberto Lanconelli\thanks{Dipartimento di Matematica, Universit\`a di Bari Aldo Moro, Bari, Italy.
\textbf{e-mail}: alberto.lanconelli@uniba.it} \and Andrea Pascucci\thanks{Dipartimento di Matematica,
Universit\`a di Bologna, Bologna, Italy. \textbf{e-mail}: andrea.pascucci@unibo.it}\and Sergio
Polidoro\thanks{Dipartimento di Scienze Fisiche, Informatiche e Matematiche, Universit\`a di Modena e Reggio Emilia,
Modena, Italy. \textbf{e-mail}: sergio.polidoro@unimore.it}}

\date{This version: \today}

\maketitle

\begin{abstract}
We prove Gaussian upper and lower bounds for the fundamental solutions of a class of degenerate
parabolic equations satisfying a weak H\"ormander condition. The bound is independent of the
smoothness of the coefficients and generalizes classical results for uniformly parabolic
equations.
\end{abstract}

\noindent \textbf{Keywords}: Kolmogorov equations, fundamental solution, linear stochastic equations, Harnack
inequalities.

%
%

\section{Introduction}
We consider the Kolmogorov backward 
equation 
\begin{equation}\label{PDE}
 \L u:=\sum_{i,j=1}^{m_0}\partial_{x_i}(a_{ij}\partial_{x_j}u)+\sum_{i=1}^{m_0}\left(\partial_{x_i}(a_{i}
u)+b_{i}\p_{x_{i}}u\right)+  c u+\sum_{i,j=1}^{d}b_{ij}x_j\partial_{x_i}u+\partial_tu=0,
\end{equation}
where $(t,x)\in\R\times\R^{d}$, $m_0\leq d$ and $\L$ verifies the following two standing assumptions:
\begin{assumption}\label{assA}
The coefficients $a_{ij}=a_{ji},a_i,b_{i},c$, for $1\le i,j\le m_0$, are bounded, measurable functions of
$(t,x)\in\mathbb{R}\times\mathbb{R}^d$ and
\begin{equation}\label{ellipticity}
 \m^{-1}|\x|^{2}\le \sum_{i,j=1}^{m_{0}}a_{ij}(t,x)\x_{i}\x_{j}\le \m|\x|^{2},\qquad
 \x\in\R^{m_{0}},\ (t,x)\in\R^{d+1},
\end{equation}
for some positive constant $\m$.
\end{assumption}
\begin{assumption}\label{assB}
The matrix $B:=\left(b_{ij}\right)_{1\leq i,j\leq d}$ has constant real entries and takes the
block-form
\begin{equation}\label{e65b}
  B=\begin{pmatrix}
 \ast & \ast & \cdots & \ast & \ast \\ B_1 & \ast &\cdots& \ast & \ast \\ 0 & B_2 &\cdots& \ast& \ast \\ \vdots & \vdots
 &\ddots& \vdots&\vdots \\ 0 & 0 &\cdots& B_{\nu}& \ast
  \end{pmatrix}
\end{equation}
where each $B_i$ is a $\left(m_{i}\times m_{i-1}\right)$-matrix of rank $m_{i}$ with
\begin{equation}
 m_0\geq m_1\geq \cdots \geq m_{\nu}\geq 1, \qquad \sum_{i=0}^{\nu} m_i = d,
\end{equation}
and the blocks denoted by {\rm ``$\ast$''} are arbitrary.
\end{assumption}

Our main result extends the bounds proved in \cite{Aronson} and \cite{Moser3, Moser3bis} for
uniformly parabolic operators with measurable coefficients: we refer to \cite{Fabes1993} for a
description of the development of this theory for non-degenerate parabolic operators, which
includes the fundamental contributions in \cite{Nash} and \cite{Davies}.
\begin{theorem} \label{th-main}
Let $\L$ be an operator in the form \eqref{PDE}, satisfying Assumptions \ref{assA} and \ref{assB}.
Let $I = ]T_0, T_1[$ be a bounded interval. Then,  there exist four positive constants $\lambda^+,
\lambda^-, C^+, C^-$ such that
\begin{equation} \label{eq-bounds}
 C^- \Gamma^{\lambda^-}(t,x;T,y) \le \Gamma(t,x;T,y) \le C^+ \Gamma^{\lambda^+} (t,x;T,y)
\end{equation}
for every $(t,x), (T,y) \in \R^{d+1}$ with $T_0 < t < T < T_1$. The constants $\lambda^+,
\lambda^-, C^+, C^-$ depend only on $d$, $\L$ and $T_1-T_0$. In \eqref{eq-bounds}
$\Gamma^{\lambda^-}$ and $\Gamma^{\lambda^+}$ denote the fundamental solutions of $\L^{\lambda^-}$
and $\L^{\lambda^+}$, respectively, where
\begin{equation} \label{eq-L-lambda}
 \L^{\lambda}:= \tfrac{\lambda}{2}\sum_{i=1}^{m_{0}}\p_{x_{i}x_{i}} + \sum_{i,j=1}^{d}b_{ij}x_j\partial_{x_i} +
 \partial_t.
\end{equation}
The explicit expression of $\Gamma^{\lambda^{\pm}}$ is given in \eqref{e22-L-lambda} below.
\end{theorem}

Degenerate equations of the form \eqref{PDE} naturally arise in the theory of stochastic processes, in physics and in
mathematical finance. For instance, if $W$ denotes a real Brownian motion, then the simplest non-trivial Kolmogorov
operator
  $$\tfrac{1}{2} \p_{vv}+v\p_{x}+\p_{t},\qquad t\ge 0,\,(v,x)\in\R^{2},$$
is the infinitesimal generator of the classical Langevin's 
stochastic equation
  $$
  \begin{cases}
    dV_{t}=dW_{t}, \\
    dX_{t}=V_{t}dt,
  \end{cases}
  $$
that describes the position $X$ and velocity $V$ of a particle in the phase space (cf. \cite{Langevin}). Notice that in
this case we have $1=m_{0}<d=2$.

Linear Fokker-Planck equations (cf. \cite{Desvillettes} and \cite{Risken}), non-linear Boltzmann-Landau equations (cf.
\cite{Lions1} and \cite{Cercignani}) and non-linear equations for Lagrangian stochastic models commonly used in the
simulation of turbulent flows (cf. \cite{Talay}) can be written in the form
\begin{equation}\label{PDE1}
 \sum_{i,j=1}^{n}\partial_{v_i}(a_{ij}\partial_{v_j}f)+\sum_{j=1}^{n}v_{j}\p_{x_{j}}f+\p_{t}f=0,\qquad
 t\ge 0,\, v\in\R^{n},\, x\in\R^{n},
\end{equation}
with the coefficients $a_{ij}=a_{ij}(t,v,x,f)$ that may depend on the solution $f$ through some integral expressions.
Clearly \eqref{PDE1} is a particular case of \eqref{PDE} with $n=m_{0}<d=2n$ and
  $$B=\begin{pmatrix}
    0_{n} & 0_{n} \\
    I_{n} & 0_{n} \
  \end{pmatrix}$$
where $I_{n}$ and $0_{n}$ denote the $\left(n\times n\right)$-identity matrix and the $\left(n\times
n\right)$-zero matrix, respectively.

In mathematical finance, equations of the form \eqref{PDE} appear in various models for the pricing of path-dependent
derivatives such as Asian options (cf., for instance, \cite{pascuccibook}, \cite{BarucciPolidoroVespri}), stochastic
volatility models (cf. \cite{HobsonRogers}, \cite{Peszek}) and in the theory of stochastic utility (cf.
\cite{AntonelliBarucciMancino}, \cite{AntonelliPascucci}). 
%
\medskip

Besides its applicative interest, the operator $\L$ in \eqref{PDE} has been studied by
several authors because of its challenging theoretical features. As in the study of uniformly parabolic operators, the theoretical results mainly depend on the assumptions on the coefficients. 
We summarize here the main results available in the literature and we focus in particular on those that are useful for the purpose of this work:
\begin{itemize}
  \item[-] {\it Constant coefficients.} If the $a_{ij}$'s, the $a_i$'s and the $b_i$'s are constant and $c=0$, the
operator $\L$ appears as the prototype of \emph{hypoelliptic operators} in the seminal
H\"ormander's work \cite{Hormander}. In particular, H\"ormander proves that a smooth fundamental
solution for $\L$ exists if, and only if, Assumptions \ref{assA} and \ref{assB} are satisfied. We
emphasize that this regularity property is not obvious for strongly degenerate operators of the
form \eqref{PDE}. Based on the explicit expression of the fundamental solution, mean value
formulas and Harnack inequalities for the non-negative solutions of $\L u = 0$ have been proved in
\cite{Kupcov5, Kupcov4, GarofaloLanconelli, LanconelliPolidoro}. In particular,
\cite{LanconelliPolidoro} studies the invariance of the solutions of $\L u = 0$ with respect to
suitable \emph{non-Euclidean} translations and \emph{non-homogeneous} dilations: it is then proved
a Harnack inequality which is translation- and dilation-invariant. In Section \ref{sec-2} we give
the precise statement of the above assertions.
  \item[-] {\it H\"older continuous coefficients.} The existence of a fundamental solution for operators
$\L$ with H\"older continuous coefficients has been proved by several authors using the parametrix
method \cite{Weber, Il'in, Sonin, Polidoro2, DiFrancescoPascucci2, DelarueMenozzi}. An invariant
Harnack inequality has been proved in \cite{Polidoro2, DiFrancescoPolidoro} and a lower bound for
the fundamental solution of $\L$ is obtained in \cite{Polidoro1, DiFrancescoPolidoro}.
  \item[-] {\it Measurable coefficients.} An upper bound for the fundamental solution of $\L$ is obtained in
\cite{PP2003, LanconelliPascucci2016} by adapting the Aronson's method \cite{Aronson}. It is based
on a local $L^\infty$-estimate of the solutions based on a Moser's iterative procedure which in
turn relies on the combination of a Caccioppoli inequality with a Sobolev estimate (see
\cite{PP2004, CPP2008, LanconelliPascucci2016}). The authors of \cite{WangZhang2011} prove a weak
form of the Poincar\'e inequality which yields the $C^\alpha$-regularity of the solutions of $\L u
= 0$. More recently, \cite{GolseVasseur} and \cite{ImbertMouhot} independently provide an
alternative proof of the $C^\alpha$-regularity of the weak solutions. Later, in the joint work
\cite{4authors}, the aforementioned four authors use their result to prove an invariant Harnack
inequality for the positive  solutions of  \eqref{PDE1}.
\end{itemize}

The main result of this paper is a lower bound for the fundamental solution $\Gamma$ of $\L$ by
merely assuming the measurability and boundedness of its coefficients, in the spirit of the works
\cite{Aronson} and \cite{Moser3, Moser3bis}. Its proof is based on the repeated application of the
Harnack inequality on suitable sequences of points that are usually called \emph{Harnack chains}.
The method used in \cite{4authors} seems us to be appropriate for the more general class of
operators $\L$ in \eqref{PDE} satisfying Assumptions \ref{assA} and \ref{assB}; for this reason,
we assume the validity of the invariant Harnack inequality for the positive weak solutions of $\L
u = 0$ and use it to prove the lower bound in \eqref{eq-bounds} under these assumptions. This
choice allows us to point out more clearly the geometric aspect of our argument. We also mention
the forthcoming paper \cite{EleuteriPolidoro} which aims at extending the techniques of
\cite{4authors} to the more general class of the operators satisfying Assumptions \ref{assA} and
\ref{assB}.

%

\section{Preliminaries} \label{sec-2}
Hereafter the operator $\L$ in \eqref{PDE} will be written in the compact form
\begin{equation}\label{L}
 \L u=\div(AD  u+au)+\langle b,Du\rangle+cu+Yu=0,
\end{equation}
{where $D=(\p_{x_{1}},\dots,\p_{x_{d}})$ denotes the gradient in $\R^{d}$,}
$A:=\left(a_{ij}\right)_{1\leq i,j\leq d}$, $a:=\left(a_i\right)_{1\leq i\leq d}$,
$b:=\left(b_i\right)_{1\leq i\leq d}$ with $a_{ij}=a_{i}=b_{i}\equiv 0$ for $i>m_0$ or $j>m_0$,
and
\begin{equation}\label{def Y}
 Y:=
 \langle B x, D   \rangle + \partial_t. 
\end{equation}
The {\it constant-coefficient Kolmogorov operator}
\begin{equation}\label{e14b}
 \LL:= \tfrac{1}{2}\sum_{i=1}^{m_{0}}\p_{x_{i}x_{i}}+ Y
\end{equation}
will be referred to as the \emph{principal part of $\L$}. It will be clear in the sequel that $\LL$ plays in this
setting the role played by the heat operator in the uniformly parabolic case. We focus here, in particular, on the
regularity properties of $\LL$ and on its invariance with respect to a family of non-Euclidean translations and
non-homogeneous dilations. It is known that Assumption \ref{assB} is equivalent to the hypoellipticity of $\LL$; in
fact, Assumption \ref{assB} is also equivalent to the well-known H\"ormander's condition, which in our setting reads:
\begin{align}\label{ae1}
 {\rm rank\ Lie}\left(\p_{x_{1}},\dots,\p_{x_{m_{0}}},Y\right)(t,x)=d+1,\qquad\text{for all }(t,x)\in \R^{d+1},
\end{align}
where Lie$\left(\p_{x_{1}},\dots,\p_{x_{m_{0}}},Y\right)$ denotes the Lie algebra generated by the
vector fields $\p_{x_{1}},\dots,\p_{x_{m_{0}}}$ and $Y$ (see Proposition 2.1 in
\cite{LanconelliPolidoro}). Thus operator $\L$ can be regarded as a perturbation of its principal
part $\LL$: roughly speaking, Assumption \ref{assA} ensures that the sub-elliptic structure of
$\LL$ is preserved under perturbation.

Constant-coefficient Kolmogorov operators are naturally associated to {\it linear} stochastic
differential equations: indeed, $\LL$ is the infinitesimal generator of the $d$-dimensional SDE
\begin{equation}\label{SDE}
  dX_{t}=B X_{t}dt+\s dW_{t},
\end{equation}
where $W$ is a standard $m_{0}$-dimensional Brownian motion and $\s$ is the
$(d\times m_{0})$-matrix 
\begin{equation}\label{sigm}
 \s=\begin{pmatrix}
   I_{m_{0}} \\
   0 \
 \end{pmatrix}.
\end{equation}
The solution $X$ of \eqref{SDE} is a Gaussian process with transition density
\begin{align} \label{e22and}
 \Gg(t,x;T,y)
 &=  \frac{1}{  \sqrt{(2\pi)^{d}\det \mathcal{C}(T-t)} }
    \exp\left(-\frac{1}{2}\langle\mathcal{C}(T-t)^{-1} \big(y - e^{(T-t)B}x)\big),
    \big(y -e^{(T-t)B}x\big)\rangle\right)
\end{align}
for  $t<T$ and $x,y\in\R^{d}$, where
\begin{align}\label{ae5}
 \mathcal{C}(t)=\int\limits_{0}^{t}\left(e^{s B}\s\right)\left(e^{s B}\s\right)^{\ast}ds
\end{align}
is the covariance matrix of $X_{t}$. Assumption \ref{assB} ensures (actually, is equivalent to the
fact) that $\mathcal{C}(t)$ is positive definite for any positive $t$. Moreover $\Gg$ in
\eqref{e22and} is the fundamental solution of $\LL$ and the function
\begin{align}
 u(t,x):=E\left[\phi\left(X_{T}\right)\mid X_{t}=x\right]=\int_{\R^{d}}\Gg(t,x;T,y)\phi(y)dy,\qquad t<T,\ x\in\R^{d},
\end{align}
solves the backward Cauchy problem
\begin{equation}\label{PC}
  \begin{cases}
    \LL u(t,x)=0,\qquad  & t<T,\ x\in\R^{d}, \\
    u(T,x)=\phi(x) & x\in\R^{d},
  \end{cases}
\end{equation}
for any bounded and continuous function $\phi$.

Operator $\LL$ has some remarkable invariance properties that were first studied in
\cite{LanconelliPolidoro}. Denote by $\ell_{(\t,\x)}$, for $(\t,\x)\in\R^{d+1}$, the
left-translations in $\R^{d+1}$ defined as
\begin{equation}\label{eq:translation}
 \ell_{(\t,\x)}(t,x):=(\t,\x)\circ (t,x):=\left(t+\t,x+{e^{t B}}\x\right),
\end{equation}
Then, $\LL$ is invariant with respect to $\ell_{\z}$ in the sense that
 $$\LL \left(u\circ \ell_{\z}\right)=\left(\LL u\right)\circ \ell_{\z},\qquad \z\in\R^{d+1}.$$
Moreover, let $\Dil(r)$ be defined as
\begin{equation}\label{e15}
 {\Dil(r) :=\text{\rm diag}(r I_{m_{0}},r^{3} I_{m_{1}},\dots,r^{2{\nu}+1}I_{m_{{\nu}}}),\qquad r\ge 0,}
\end{equation}
where $I_{m_{i}}$ denotes the $\left(m_{i}\times m_{i}\right)$-identity matrix. Then, $\LL$ is homogeneous with respect
to the dilations in $\R^{d+1}$ defined as
\begin{equation}\label{e066}
 \d_{r}(t,x):=\left(r^2t,\Dil(r)x\right),
\end{equation}
{\it if and only if} all the $\ast$-blocks of $B$ in \eqref{e65b} are null
(\cite{LanconelliPolidoro}, Proposition 2.2). In this case, we have
\begin{equation}\label{e200}
 \LL(u\circ\delta_r)=r^{2}(\LL u)\circ\delta_r.
\end{equation}
The natural number
\begin{equation}\label{e77}
 Q:=m_{0}+3m_{1}+\cdots+(2\nu+1)m_{\nu}.
\end{equation}
is usually called the {\it homogeneous dimension} of $\mathbb{R}^{d}$ with respect to
$(\Dil(r))_{r>0}$, because the Jacobian of $\Dil(r)$ is equal to $r^{Q}$.

In accordance with \eqref{e22and}, the fundamental solution of the operator $\L^{\lambda}$ defined in
\eqref{eq-L-lambda} is
\begin{align} \label{e22-L-lambda}
 \Gamma^{\lambda}(t,x;T,y)
 &=  \frac{1}{  \sqrt{(2\pi \lambda)^{d}\det \mathcal{C}(T-t)} }
    \exp\left(-\frac{1}{2 \lambda}\langle\mathcal{C}(T-t)^{-1} \big(y - e^{(T-t)B}x)\big),
    \big(y -e^{(T-t)B}x \rangle\right)
\end{align}
for  $t<T$ and $x,y\in\R^{d}$.

We end this section with the definitions of weak and fundamental solutions utilized in the sequel.

\begin{definition}
A \emph{weak solution} of \eqref{PDE} in a domain $\O$ of $\R^{d+1}$ is a function $u$ such that
  $$u,\p_{x_{1}}u,\dots, \p_{x_{m_{0}}}u, Yu\in L^{2}_{\text{\rm loc}}(\Omega)$$
and
\begin{equation}
 \int_{\Omega}-\langle AD  u,D  \psi\rangle-u\langle a,D  \psi\rangle +\psi \langle b,Du\rangle+u\psi+\psi
 Yu=0,
\end{equation}
for any $\psi\in C_0^{\infty}(\Omega)$.
\end{definition}


We recall that the formal adjoint operator of $\L$ is defined as
\begin{equation}\label{PDEagg}
 \L^*v:=\sum_{i,j=1}^{m_0}\partial_{y_i}(a_{ij}\partial_{y_j}v)-\sum_{i=1}^{m_0}\left(a_{i}\partial_{y_i}
v+\p_{y_{i}}(b_{i}v)\right)+
 \left(c- \text{tr}(B)\right) v-\sum_{i,j=1}^{d}b_{ij}y_j\partial_{y_i}v-\partial_T v.
\end{equation}
\begin{definition}
A fundamental solution for $\L$ is a continuous and positive function $\G=\G(t,x;T,y)$, defined for
$t<T$ and $x,y\in\R^{d}$, such that:
\begin{itemize}
  \item[i)] $\G( \cdot,\cdot;T,y)$ is a weak solution of $\L u=0$ in $]-\infty,T[\times \R^{d}$ and
  $\G(t,x;\cdot,\cdot)$ is a weak solution of $\L^{*} u=0$ in $]t,+\infty[\times \R^{d}$;
  \item[ii)] for any bounded function $\phi\in C(\R^{d})$ and $x,y\in\R^{d}$, we have
\begin{align}
  \lim_{(t,x)\to(T,y)\atop t<T}u(t,x)=\phi(y), 
  \qquad \lim_{(T,y)\to(t,x)\atop T>t}v(T,y)=\phi(x), 
\end{align}
where
  \begin{align}\label{ae11}
 u(t,x):=\int_{\R^{d}}\G(t,x;T,y)\phi(y)dy,\qquad 
 v(T,y):=\int_{\R^{d}}\G(t,x;T,y)\phi(x)dx.
\end{align}
\end{itemize}
\end{definition}
\begin{remark}
The functions in \eqref{ae11} are weak solutions of the following backward and forward Cauchy problems:
  $$
  \begin{cases}
    \L u(t,x)= 0,\ &(t,x)\in\, ]-\infty,T[\times \R^{d}, \\
    u(T,x)= \phi(x), & x \in\R^{d},
  \end{cases}\qquad
  \begin{cases}
    \L^{*}v(T,y)= 0,\ & (T,y)\in\,]t,+\infty[\times \R^{d}, \\
    v(t,y)= \phi(y) & y \in\R^{d}.
  \end{cases}
  $$
\end{remark}

\begin{remark}
A fundamental solution for $\L$ exists under the additional condition that the coefficients are
H\"older continuous (see \cite{Polidoro2, DiFrancescoPascucci2, DelarueMenozzi}) .
\end{remark}

\begin{remark}\label{l3}
Let $u$ be a weak solution of \eqref{PDE} and $r>0$. Then $v:=u\circ\d_{r}$ solves $\L^{(r)}v=0$
where
\begin{equation}\label{Lr}
  \L^{(r)}v:=\div(A^{(r)}D v)+\div(a^{(r)}v)+\langle b^{(r)},D v\rangle+c^{(r)}v+\langle B^{(r)}x,D
v\rangle+\partial_t v,
\end{equation}
with $A^{(r)}=A\circ \d_{r}$, $a^{(r)}=r (a\circ \d_{r})$, $b^{(r)}=r (b\circ \d_{r})$,
$c^{(r)}=r^{2}(c\circ \d_{r})$
and $B^{(r)}=r^{2}\Dil_{r} B\Dil_{\frac{1}{r}}$, that is
\begin{equation}\label{e65bl}
  B^{(r)}=\begin{pmatrix}
 r^{2} B_{1,1} & r^{4} B_{1,2} & \cdots & r^{2\nu} B_{1,\nu} & r^{2\nu+2} B_{1,\nu+1} \\
 B_1 & r^{2} B_{2,2} &\cdots& r^{2\nu-2} B_{2,\nu} & r^{2\nu} B_{2,\nu+1} \\
 0 & B_2 &\cdots& r^{2\nu-4} B_{3,\nu}& r^{2\nu-2} B_{3,\nu+1} \\ \vdots & \vdots
 &\ddots& \vdots&\vdots \\ 0 & 0 &\cdots& B_{\nu}& r^{2} B_{\nu+1,\nu+1}
  \end{pmatrix},
\end{equation}
where $B_{i,j}$ denotes the $\ast$-block in the $(i,j)$-th position in \eqref{e65b}.
\end{remark}

\begin{notation}\label{notation1} Let $\costd>0$ and $B:=\left(b_{ij}\right)_{1\leq i,j\leq
d}$ a matrix that satisfies Assumption \ref{assB}. We denote by $\Kol$ the class of Kolmogorov
operators of the form \eqref{Lr} with $\l\in[0,1]$ and the coefficients $a_{ij},a_i,b_{i},c$, for
$1\le i,j\le m_0$, that satisfy Assumption \ref{assA} with the non-degeneracy constant $\m$ in
\eqref{ellipticity} and the norms
$\|a_{i}\|_{\infty}$, $\|b_{i}\|_{\infty}$, $\|c\|_{\infty}$ smaller than $\costd$. 
\end{notation}
\begin{remark}\label{r7}
Let $\L \in\Kol$. If $u$ is a solution of $\L u=0$ then, for any $\z\in\R^{d+1}$, $v:=u\circ
\ell_{\z}$ solves  $\left(\L \circ \ell_{\z}\right)v=0$ where $(\L \circ \ell_{\z})$ is the operator
obtained from $\L$ by $\ell_{\z}$-translating its coefficients. Moreover, operator
$(\L \circ \ell_{\z})$ still belongs to $\Kol$.
\end{remark}


\section{Harnack inequalities}
Let $B$ be a matrix that satisfies Assumption \ref{assB}. We associate to $B$ the cylinders
\begin{align}
 Q^+_{1} = \{(t,x)\in\R\times\R^{d}\mid  0\le t<1,\, |x|<1\},
\end{align}
and
\begin{equation}\label{cylinder}
 Q^+_{r}(z_{0}):= z_0 \circ \d_{r} \left(Q^+_{1}\right)=\{z\in\R^{d+1}\mid z= z_{0} \circ \d_{r}(\z), \, \z \in
 Q^+_{1}\},
\end{equation}
for $z_{0}\in\R^{d+1}$ and $r>0$. The following remarkable result is proved in
\cite{4authors} for prototype Kolmogorov equations \eqref{PDE1} and in \cite{EleuteriPolidoro} for
general Kolmogorov equations  \eqref{PDE}.
\begin{theorem}[{\bf Local Harnack inequality}]\label{t-h1}
Let $\L \in \Kol$. If $u$ is a non-negative weak solution of \eqref{PDE} in $Q^+_1$ then
\begin{equation}\label{h4a}
 \sup_{Q^+_r(\beta,0)}u \leq C\inf_{Q^+_r(0,0)}u,
\end{equation}
where the constants $C\ge 1$ and $\b,r\in\,]0,1[$ depend only on $M$ and $B$.
\end{theorem}
\begin{remark}
The constants $\b,r$ in Theorem \ref{t-h1} are small so that the cylinders $Q^+_r(0,0)$ and
$Q^+_r(\b,0)$ are disjoint subsets of $Q^+_1$.
\end{remark}
\begin{remark}\label{r5}
By Remark \ref{r7}, the Harnack inequality \eqref{h4a} is valid for cylinders centered at an
arbitrary point $z_{0}\in\R^{d}$ with the same constants $C,\b,r$, dependent only on $M$ and $B$.
\end{remark}

Next we prove a global version of the Harnack inequality based on a classical argument which makes
use of the so-called Harnack chains. We first prove a preliminary result. For $\b,r,R>0$ and
$z_{0}\in\R^{d+1}$, we define the cones
\begin{equation}\label{e-cone}
    \P_{\b,r,R}  = \left\{ z\in\R^{d+1}\mid z=\d_{\l}(\b,\x),\, |\x| < r,\ 0<\l\le R \right\},
\end{equation}
and $\P_{\b,r,R} (z_{0}) := z_{0} \circ \P_{\b,r,R}$. Here $ |\x|$ denotes the Euclidean norm of the vector
$\x \in \R^d$. Theorem \ref{t-h1} combined with Remark \ref{l3} gives the following
\begin{lemma}\label{c6}
Let $z\in\R^{d+1}$ and $R\in\,]0,1]$. Let $u$ be a continuous and non-negative weak solution of \eqref{PDE} in
$Q^+_{R}(z)$. Then we have
\begin{equation}\label{h4abis}
 \sup_{P_{\b,r,R}(z)}u \leq C u(z),
\end{equation}
where the constants $C,\b$ and $r$ are the same as in Theorem \ref{t-h1} and depend only on $M$
and $B$.
\end{lemma}

\proof We preliminarily recall that \cite{WangZhang2011} prove the H\"older continuity of the solutions of $\L u = 0$.
In particular, $u$ is a continuous function. So let $u$ be a continuous and non-negative weak solution of \eqref{PDE}
in $Q^+_{R}(z)$ and let $w \in P_{\b,r,R}(z)$. Then $w=z \circ \d_{\l}(\b,\x)$ for some $\l\in\,]0, R]$ and $|\x| < r$.
By using the notation introduced in \eqref{eq:translation}, we obtain from Remark \ref{l3} that the function $u_{z,
\l}:=u\circ  \ell_{z} \circ \d_{\l}$ is a continuous and non-negative weak solution in $Q^+_{\frac{R}{\l}}(0,0)
\supseteq Q^+_{1}(0,0)$ of $\L^{(\l)} u_{\l} = 0$, where $\L^{(\l)}$ is the operator defined in \eqref{Lr}. Since
$\L^{(\l)}\in\Kol$, by the Harnack inequality \eqref{h4a} for $\L^{(\l)}$, we have
\begin{align}
 u(w) = u_{z,\l}(\b,\x) \le \sup_{Q^{+}_{r}(\b,0)} u_{z,\l}\le C \inf_{Q^{+}_{r}(0,0)} u_{z,\l}\le
  C u_{z,\l}(0,0)=C u(z).
\end{align}
\endproof

%
%

\begin{theorem}[{\bf Global Harnack inequality}]\label{t-h3}
Let $\L \in \Kol$, $T\in\R$ and $\t\in\,]0,1]$. If $u$ is a non-negative weak solution of
\eqref{PDE} in $]T -\t,T + \t[\times\rd$, then we have
 $$u(T,y) \leq \costc e^{\costc \langle \mathcal{C}^{-1}(T-t) (y- e^{(T-t)B}x) ,
 y- e^{(T-t)B}x \rangle} u(t,x),\qquad t\in\,]T-\t, T[,\ x,y\in\R^{d},
 $$
where $\mathcal{C}$ is the covariance matrix in \eqref{ae5} and  $\costc$ is a positive constant
that depends only on $M$ and $B$.
\end{theorem}

Before proving Theorem \ref{t-h3}, we recall (see, for instance, Sect.9.5 in \cite{pascuccibook}) that the H\"ormander
condition \eqref{ae1} 
is equivalent to the fact that the pair of matrices $(B,\s)$, with $\s$ as in \eqref{sigm}, is controllable in
the following sense: for any
$(t,x), (T,y) \in \rdd$ with $t<T$, there exists $v\in L^{2}([t,T];\R^{m_{0}})$ 
such that the system
\begin{equation}\label{e-gamma2}
    \begin{cases}
    \g'(s) = B\g(s) + \sigma v(s), \\
    \g(t) = x, \quad \g(T) = y,
    \end{cases}
\end{equation}
has solution. The function $v$ is called a {\it control for $(B,\s)$ on $[t,T]$}. In the proof of
Theorem \ref{t-h3} we will use the following
\begin{lemma}\label{p-h4} Let $\g$ be the solution of the linear problem
\begin{equation}\label{e-gamma2bis}
    \begin{cases}
    \g'(s) = B\g(s) + \sigma v(s),\qquad s\in[t,T], \\
    \g(t) = x,
    \end{cases}
\end{equation}
with $T-t\le 1$, initial datum $x\in\R^{d}$ and control function $v\in L^{2}([t,T];\R^{m_{0}})$.
Then we have
 $$(s, \g(s)) \in \P_{1,\kappa \|v\|_{L^{2}([t,T])},\sqrt{T-t}} (t,x),\qquad s \in [t,T],$$
where $\kappa$ is a positive constant which depends only on $B$.
\end{lemma}
\proof The explicit solution of \eqref{e-gamma2bis} is
\begin{equation}\label{e-explicit}
    \g(s) = e^{(s-t)B}x + \int_t^s e^{(s-\tau)B} \sigma v(\tau) d \tau,\qquad s \in [t,T].
\end{equation}
Thus, setting $\l=\sqrt{s-t}$, 
we have that $(s,\g(s))\in \P_{1,r,\sqrt{T-t}}(t,x)$ if and only if
\begin{equation}\label{ae3}
  \int_t^{t+\l^{2}} e^{(t+\l^{2}-\t)B} \sigma v(\t) d \t=\mathcal{D}\left(\l\right)\x\quad
  \text{ with }\quad |\x|\le r.
\end{equation}
To check this, we first notice that, according to \eqref{e15}, the space $\R^{d}$ admits a natural
decomposition as a direct sum
  $$\R^{d}=\bigoplus_{j=0}^{\n}V_{j},\qquad \text{dim}\,V_{j}=m_{j}.$$
Then, for $x\in\R^{d}$, with obvious notation we have $x=x^{(0)}\oplus\cdots\oplus x^{(\v)}$ where
  $$\Dil(r)x^{(j)}=r^{2j+1}x^{(j)},\qquad j=0,\dots,\n.$$
We also write a $(d\times d)$-matrix $E$ in block form as in \eqref{e65b}, that is
$E=\left(E^{(ij)}\right)_{i,j=0\dots,\n}$ where $E^{(ij)}$ is a block of dimension $m_{i}\times
m_{j}$. In particular, given the definition of exponential $E(t):=e^{t B}$ as the sum of a power
series, a direct computation shows that
\begin{equation}\label{ae2}
\begin{split}
  E^{(00)}(t)&=I_{m_{0}}+t\text{O}(t),\\
  E^{(0j)}(t)&=\frac{t^{j}}{j!}\left(I_{m_{j}}+t\text{O}(t)\right)B_{j}\cdots B_{1},\qquad
  j=1,\dots,\n,
\end{split}
\end{equation}
as $t\to 0$, where $I_{m_{j}}$ denotes the $(m_{j}\times m_{j})$-identity matrix. Now, $\sigma
v\in V_{0}$ and therefore, by \eqref{ae2}, we have
  $$\left|\left(e^{(t+\l^{2}-\t )B} \sigma v(\t )\right)^{(j)}\right|\le \kappa (t+\l^{2}-\t )^{j} |v(\t )|,\qquad \t
\in[t,T],$$
with the constant $\kappa$ dependent only on $B$. Thus we have
\begin{align}
 \left|\int_t^{t+\l^{2}} \left(e^{(t+\l^{2}-\t )B} \sigma v(\t )\right)^{(j)}\right| d \t &\le
  \kappa \int_t^{t+\l^{2}} (t+\l^{2}-\t )^{j} |v(\t )| d \t \le
\intertext{(by H\"older's inequality)}
  &\le \kappa \|v\|_{L^{2}([t,T])}\l^{2j+1},
\end{align}
and this proves \eqref{ae3}.
\endproof

Let us consider the control problem \eqref{e-gamma2} one more time. Among the paths $\g$
satisfying \eqref{e-gamma2}, one is often interested in one minimizing the {\it total cost}
\begin{equation}\label{e-cost}
    \costv = \int_t^{T} |v(s)|^2 d s.
\end{equation}
Classical control theory provides the explicit expression of an optimal control and of its cost
(see, for instance, \cite{pascuccibook}, Theor. 9.55). 

\begin{lemma}\label{la1}
 The optimal control for problem \eqref{e-gamma2} is given by
\begin{equation}
     \bar{v}(s) = \left(e^{-(s-t) B}\sigma\right)^{T}  \mathcal{C}^{-1}(-(T-t)) \left(e^{-(T-t)B} y -
     x\right),\qquad s\in[t,T].
\end{equation}
The corresponding minimal cost will be denoted by
\begin{equation*}
 V(t,x;T,y) := \|\bar{v}\|_{L^{2}([t,T])}^{2}
\end{equation*}
and is equal to
\begin{equation*}
 V(t,x;T,y) = \langle \mathcal{C}^{-1}(T-t) (y- e^{(T-t)B} x) , y- e^{(T-t)B}x \rangle.
\end{equation*}
\end{lemma}
\proof[\noindent Proof of Theorem \ref{t-h3}] In order to use the previous versions of the Harnack
inequality, we first notice that by assumption,  for every $z\in ]T-\t,T[\times\R^{d}$, $u$ is a
continuous and non-negative weak solution of \eqref{PDE} in $Q^{+}_{\sqrt{\t}}(z)$. Next we fix $x,y\in\R^{d}$,
$t\in\,]T-\t,T[$ and consider the solution $\g$ of the control problem \eqref{e-gamma2} corresponding to the optimal
control $\bar{v}$ given in Lemma \ref{la1}. Moreover, we set $\coste=\left(\frac{r}{\kappa}\right)^{2}$ where $r$ and
$\kappa$ are the constants in Theorem \ref{t-h1} and Lemma \ref{p-h4} respectively.

Now, if $T\le t+\t\b$ and $\|\bar{v}\|^{2}_{L^{2}([t,T])}\le \coste$, then by Lemma \ref{p-h4} we have
  $$(T,y)\in P_{1,r,\sqrt{\t}}(t,x)\cap\left(]t,t+\t\b]\times\R^{d}\right)\subseteq P_{\b,r,\sqrt{\t}}(t,x)$$
and therefore by Lemma \ref{c6} we get
  $$u(T,y)\le C u(t,x)$$
where $C$ is the constant in Theorem \ref{t-h1}, which depends only on $M$ and $B$.

Viceversa, setting $t_{0}=t$ and
  $$t_{j+1}=(t_{j}+\t\b)\wedge\inf\{s\in[t_{j},T]\mid \|\bar{v}\|^{2}_{L^{2}([t_{j},s])}
  \ge \coste\},$$
we have that $t_{j}=T$ for $j\ge \frac{1}{\b}+\frac{\|\bar{v}\|^{2}_{L^{2}([t,T])}}{\coste}$ and
  $$(t_{j+1},\g(t_{j+1})) \in P_{1,r,\sqrt{\t}}(t_{j},\g(t_{j}))\cap\left(]t_{j},t_{j}+\t\b]\times\R^{d}\right)
  \subseteq P_{\b,r,\sqrt{\t}}(t_{j},\g(t_{j}))$$
if $t_{j}<T$. By Lemma \ref{c6} we have
  $$u(t_{j},\g(t_{j}))\le C u(t_{j-1},\g(t_{j-1})),$$
which yields
  $$u(T,y)\le C^{\frac{1}{\b}+\frac{1}{\coste} V(t,x;T,y) }u(t,x).$$
The thesis follows by using the expression of the optimal cost given in Lemma \ref{la1}.
\endproof

\section{Lower bounds for fundamental solutions}

The proof of the lower bound for the fundamental solution will make use of the following upper
bound obtained in \cite{PP2003} and \cite{LanconelliPascucci2016}.
\begin{theorem}[\bf Gaussian upper bound]\label{t1lp}
Let $\L\in \Kol$. There exists a positive constant $\cca$, only dependent on $\costd$ and $B$, 
such that
\begin{align}\label{thes}
 \G(t,x;T,y)\le
\frac{\cca}{(T-t)^{\frac{Q}{2}}}\exp\left(-\frac{1}{\cca}\left|\Dil\left(\left(T-t\right)^{-\frac{1}{2}}\right)
 \left(y-e^{(T-t)B}x\right)\right|^2\right),
\end{align}
for $0<T-t\le 1$ and $x,y\in\R^{d}$.
\end{theorem}

\begin{lemma}\label{lemma}
Let $\L\in \Kol$. There exist two positive constants $R$ and $\ccb$, which depend only on $\costd$
and $B$, such that
\begin{equation}\label{ae6}
 \int_{\left|\mathcal{D}(\sqrt{T-t})\left(y-e^{(T-t)B}x\right)\right|\le R}\Gamma(t,x;T,y)dx\geq \ccb,
 \qquad 0<T-t\le 1,\ y\in\R^{d}.
\end{equation}
\end{lemma}
\proof First notice that, for a suitably large constant $\ccc$ dependent only on $M$ and $B$, the
function
  $$v(T,y):=\int_{\R^{d}}\Gamma(t,x;T,y)dx-e^{-\ccc(T-t)},\qquad T>t,\ y\in\R^{d},$$
is a weak super-solution of the forward Cauchy problem
  $$
  \begin{cases}
    \L^{*}v(T,y)=-e^{-\ccc(T-t)}\left(c- \text{tr}(B)+\ccc\right)\le 0,\ & T>t,\ y\in\R^{d}, \\
    v(t,y)=0\ & y\in\R^{d},
  \end{cases}
  $$
for the adjoint operator $\L^{*}$ in \eqref{PDEagg}. Therefore, by the maximum principle (see, for
instance, Proposition 3.4 in \cite{DifrancescoPascucciPolidoro2007}), we have $v\ge 0$ that is
  $$\int_{\R^{d}}\Gamma(t,x;T,y)dx\ge e^{-\ccc(T-t)},\qquad T>t,\ y\in\R^{d}.$$
Then \eqref{ae6} follows from the following estimate:
\begin{align}
 &\int_{\left|\mathcal{D}(\sqrt{T-t})\left(y-e^{(T-t)B}x\right)\right|\ge R}\Gamma(t,x;T,y)dx\le
\intertext{(by the upper bound \eqref{thes})}
 &\le \frac{\cca}{(T-t)^{\frac{Q}{2}}}
 \int_{\left|\mathcal{D}(\sqrt{T-t})\left(y-e^{(T-t)B}x\right)\right|\ge
R}\exp\left(-\frac{1}{\cca}\left|\Dil\left(\left(T-t\right)^{-\frac{1}{2}}\right)
 \left(y-e^{(T-t)B}x\right)\right|^2\right)dx=
\intertext{(by the change of variable $z=\Dil\left(\left(T-t\right)^{-\frac{1}{2}}\right)
 \left(y-e^{(T-t)B}x\right)$)}
 &= \cca\int_{\left|z\right|\ge R}\exp\left(-\frac{1}{\cca}\left|z\right|^2\right)dz
\end{align}
which gives the thesis.
\endproof

We are now ready to state and prove the main result of the present paper.
\begin{theorem}[\bf Gaussian lower bound]\label{t1}
Let $\L\in \Kol$. There exists a positive constant $\ccd$, dependent only on $M$ and $B$, such that
\begin{equation}\label{ae7}
 \G(t,x;T,y)\geq \frac{\ccd}{(T-t)^{\frac{Q}{2}}}e^{- \frac{1}{\ccd} \langle \mathcal{C}^{-1}(T-t)
(y-e^{(T-t)}x),y-e^{(T-t)}x \rangle},\qquad 0<T-t\le 1,\ x,y\in\R^{d}.
\end{equation}
\end{theorem}
\begin{remark}\label{r8}
In general, estimate \eqref{ae7} is valid for any $T-t>0$, with $\ccd$ dependent also on
$1\vee(T-t)$.
\end{remark}
\proof We prove a preliminary diagonal estimate. Let $\t=\frac{T-t}{2}$: by the global Harnack
inequality stated in Theorem \ref{t-h3}, for any $\x,y\in\R^{d}$ we have
\begin{align}\label{ae9}
 \G(t,y;T,y)\ge \costc e^{-\costc \langle \mathcal{C}^{-1}(\t) (\x- e^{\t B}y), \x- e^{\t B}y \rangle} \G(t+\t,\x;T,y).
\end{align}
For any $y \in \R^d$ we set
  $$D_{R}=\{\x\in\R^{d}\mid \left|\mathcal{D}(\sqrt{\t})\left(y-e^{\t B}\x\right)\right|\le R\}, \qquad R>0,$$
and notice that, up to a constant dependent only on $M$ and $B$, the Lebesgue measure of $D_{R}$
equals $\t^{Q}$. We also note that, by Lemma 3.3 in \cite{LanconelliPolidoro}, $\langle
\mathcal{C}^{-1}(\t) (\x- e^{\t B}y) , \x- e^{\t B}y \rangle$ is bounded on $D_{R}$. Therefore,
integrating \eqref{ae9} over $D_{R}$, we get
\begin{align}\label{ae10}
 \G(t,y;T,y)\ge \frac{\ccf}{\t^{Q}} \int_{\left|\mathcal{D}(\sqrt{\t})\left(y-e^{\t B}\x\right)\right|\le
 R}\G(t+\t,\x;T,y)d\x\ge \frac{\ccg}{(T-t)^{\frac{Q}{2}}},
\end{align}
where the last inequality follows from Lemma \ref{lemma} and the constant $\ccg$ depends only on
$M$ and $B$. Hence, by applying again the global Harnack inequality we get
\begin{align}\label{eq-tau-2tau}
  \G(t,0;T,y)&\ge \costc e^{-\costc \langle \mathcal{C}^{-1}(\t) y,y \rangle} \G(t+\t,y;T,y)\ge
\intertext{(by \eqref{ae10})}
  &\ge \frac{\cch}{(T-t)^{\frac{Q}{2}}} e^{-\costc \langle \mathcal{C}^{-1}(\t) y,y
  \rangle}\ge \frac{\cci}{(T-t)^{\frac{Q}{2}}} e^{-\cci \langle \mathcal{C}^{-1}(T-t) y,y
  \rangle},
\end{align}
where the last inequality is a consequence of \eqref{quadratic} from Remark \ref{equivalence} below. This proves \eqref{ae7} for $x=0$; the general statement follows by the translation-invariance property of the operator $\L$.
\endproof

\begin{remark}\label{equivalence}
If we denote by $\mathcal{C}_0$ the covariance matrix appearing in the fundamental solution of the homogenous principal part of $\mathcal{L}$, then there exist $\alpha_1,...,\alpha_4,\beta_1,...,\beta_4>0$ such that for any $\tau\in ]0,1]$ and $z\in\mathbb{R}^d$
\begin{equation}\label{determinant}
\alpha_1\tau^{\mathcal{Q}}\leq\alpha_2\det(\mathcal{C}_0(\tau))\leq\det(\mathcal{C}(\tau))
\leq\alpha_3\det(\mathcal{C}_0(\tau))\leq\alpha_4\tau^{\mathcal{Q}}
\end{equation}
and
\begin{equation}\label{quadratic}
\beta_1\left|\Dil\left(\tau^{-\frac{1}{2}}\right)z\right|^2\leq\beta_2\langle\mathcal{C}_0^{-1}(\tau)z,z\rangle\leq
\langle\mathcal{C}^{-1}(\tau)z,z\rangle\leq\beta_3\langle\mathcal{C}_0^{-1}(\tau)z,z\rangle
\leq\beta_4\left|\Dil\left(\tau^{-\frac{1}{2}}\right)z\right|^2.
\end{equation}
In fact, we recall (see Proposition 2.3 in \cite{LanconelliPolidoro}) that for any $\tau>0$ one has
\begin{equation}
\mathcal{C}_0(\tau)=\Dil\left(\sqrt{\tau}\right)\mathcal{C}_0(1)\Dil\left(\sqrt{\tau}\right)
\end{equation}
and
\begin{equation}
\mathcal{C}^{-1}_0(\tau)=\Dil\left(\tau^{-\frac{1}{2}}\right)\mathcal{C}^{-1}_0(1)\Dil\left(\tau^{-\frac{1}{2}}
\right).
\end{equation}
These identities imply that
\begin{equation}
\det\mathcal{C}_0(\tau)=\det\left(\Dil\left(\sqrt{\tau}\right)
\mathcal{C}_0(1)\Dil\left(\sqrt{\tau}\right)\right)\\ =\tau^{\mathcal{Q}}\det\mathcal{C}_0(1);
\end{equation}
moreover, if $k_1$ and $k_2$ denote, respectively, the least and the greatest eigenvalue of
$\mathcal{C}^{-1}_0(1)$, we have that  $k_1 > 0$ and
\begin{equation}
k_1\left|\Dil\left(\tau^{-\frac{1}{2}}\right)z\right|^2\leq
\langle\mathcal{C}_0^{-1}(\tau)z,z\rangle \leq
k_2\left|\Dil\left(\tau^{-\frac{1}{2}}\right)z\right|^2
\end{equation}
for all $z\in\mathbb{R}^d$ and $\tau>0$. This proves the first and last inequalities in \eqref{determinant} and \eqref{quadratic}. To prove the equivalence between the matrices $\mathcal{C}^{-1}_0$ and $\mathcal{C}^{-1}$, we recall that, according to formula (3.14) in \cite{LanconelliPolidoro}, we have
\begin{equation}
\frac{\det\mathcal{C}(\tau)}{\det\mathcal{C}_0(\tau)}=1+\tau O(1),\quad\mbox{ as }\tau\to 0^+.
\end{equation}
Hence $\frac{\det\mathcal{C}(\tau)}{\det\mathcal{C}_0(\tau)}$ can be extended to a positive and
continuous function on $[0,1]$: consequently there exist two positive constants $k_3$ and $k_4$
such that
\begin{equation}\label{det}
 k_3 \det\mathcal{C}_0(\tau) \leq \det\mathcal{C}(\tau) \le k_4 \det\mathcal{C}_0(\tau), \qquad \tau\in [0,1].
\end{equation}
By the same argument we can prove that there exist two positive constants $k_5$ and $k_6$ such
that
\begin{equation}\label{cov}
k_5 \langle\mathcal{C}_0^{-1}(\tau)z,z\rangle \le \langle\mathcal{C}^{-1}(\tau)z,z\rangle \leq k_6
\langle\mathcal{C}_0^{-1}(\tau)z,z\rangle
\end{equation}
for every $z\in\mathbb{R}^d$ and $\tau\in [0,1]$ (see inequality (2.12) in
\cite{DiFrancescoPolidoro}). Indeed we recall that, for every $z \in \R^d$,
\begin{equation}
\langle\mathcal{C}^{-1}(\tau)z,z\rangle = \langle\mathcal{C}^{-1}(\tau)z,z\rangle=1+\tau O(1),\quad\mbox{ as
}\tau\to 0^+.
\end{equation}
(see Lemma 3.3 in \cite{LanconelliPolidoro}.) Then, the function $(\tau,z) \mapsto
\frac{\langle\mathcal{C}^{-1}(\tau)z,z\rangle}{\langle\mathcal{C}^{-1}(\tau)z,z\rangle}$ can be
extended to a positive and continuous function on the compact set
\begin{equation*}
  \left\{(\tau,z) \in \R^{d+1} \mid 0 \le \tau \le 1,\, |z| = 1\right\}.
\end{equation*}
Then we conclude as above.
\end{remark}

The following corollary is a straightforward consequence of Theorem \ref{t1lp} and Remark \ref{equivalence}
\begin{corollary}
Let $\L\in \Kol$. There exists a positive constant $\ccl$, only dependent on $\costd$ and $B$, such that
\begin{align}\label{thes4}
\G(t,x;T,y)\le\frac{\ccl}{\sqrt{\det\mathcal{C}(T-t)}}\exp\left(-\frac{1}{\ccl}\langle\mathcal{C}^{-1}(T-t)
\left(y-e^{(T-t)B}x\right),\left(y-e^{(T-t)B}x\right)\rangle\right).
\end{align}
for $0<T-t\le 1$ and $x,y\in\R^{d}$.
\end{corollary}

\bibliographystyle{siam}
\bibliography{BibTeX-Final}

\def\cprime{$'$} \def\cprime{$'$} \def\cprime{$'$} \def\cprime{$'$}
  \def\lfhook#1{\setbox0=\hbox{#1}{\ooalign{\hidewidth
  \lower1.5ex\hbox{'}\hidewidth\crcr\unhbox0}}} \def\cprime{$'$}
  \def\cprime{$'$} \def\cprime{$'$}
\begin{thebibliography}{10}

\bibitem{AntonelliBarucciMancino}
{\sc F.~Antonelli, E.~Barucci, and M.~E. Mancino}, {\em Asset pricing with a
  forward-backward stochastic differential utility}, Econom. Lett., 72 (2001),
  pp.~151--157.

\bibitem{AntonelliPascucci}
{\sc F.~Antonelli and A.~Pascucci}, {\em On the viscosity solutions of a
  stochastic differential utility problem}, J. Differential Equations, 186
  (2002), pp.~69--87.

\bibitem{Aronson}
{\sc D.~G. Aronson}, {\em Bounds for the fundamental solution of a parabolic
  equation}, Bull. Amer. Math. Soc., 73 (1967), pp.~890--896.

\bibitem{BarucciPolidoroVespri}
{\sc E.~Barucci, S.~Polidoro, and V.~Vespri}, {\em Some results on partial
  differential equations and {A}sian options}, Math. Models Methods Appl. Sci.,
  11 (2001), pp.~475--497.

\bibitem{Talay}
{\sc M.~Bossy, J.-F. Jabir, and D.~Talay}, {\em On conditional {M}c{K}ean
  {L}agrangian stochastic models}, Probab. Theory Related Fields, 151 (2011),
  pp.~319--351.

\bibitem{Cercignani}
{\sc C.~Cercignani}, {\em The {B}oltzmann equation and its applications},
  Springer-Verlag, New York, 1988.

\bibitem{CPP2008}
{\sc C.~Cinti, A.~Pascucci, and S.~Polidoro}, {\em Pointwise estimates for a
  class of non-homogeneous {K}olmogorov equations}, Math. Ann., 340 (2008),
  pp.~237--264.

\bibitem{Davies}
{\sc E.~B. Davies}, {\em Explicit constants for {G}aussian upper bounds on heat
  kernels}, Amer. J. Math., 109 (1987), pp.~319--333.

\bibitem{DelarueMenozzi}
{\sc F.~Delarue and S.~Menozzi}, {\em Density estimates for a random noise
  propagating through a chain of differential equations}, J. Funct. Anal., 259
  (2010), pp.~1577--1630.

\bibitem{Desvillettes}
{\sc L.~Desvillettes and C.~Villani}, {\em On the trend to global equilibrium
  in spatially inhomogeneous entropy-dissipating systems: the linear
  {F}okker-{P}lanck equation}, Comm. Pure Appl. Math., 54 (2001), pp.~1--42.

\bibitem{DiFrancescoPascucci2}
{\sc M.~Di~Francesco and A.~Pascucci}, {\em On a class of degenerate parabolic
  equations of {K}olmogorov type}, AMRX Appl. Math. Res. Express, 3 (2005),
  pp.~77--116.

\bibitem{DifrancescoPascucciPolidoro2007}
{\sc M.~Di~Francesco, A.~Pascucci, and S.~Polidoro}, {\em The obstacle problem
  for a class of hypoelliptic ultraparabolic equations}, Proc. R. Soc. Lond.
  Ser. A Math. Phys. Eng. Sci., 464 (2008), pp.~155--176.

\bibitem{DiFrancescoPolidoro}
{\sc M.~Di~Francesco and S.~Polidoro}, {\em Schauder estimates, {H}arnack
  inequality and {G}aussian lower bound for {K}olmogorov-type operators in
  non-divergence form}, Adv. Differential Equations, 11 (2006), pp.~1261--1320.

\bibitem{EleuteriPolidoro}
{\sc M.~Eleuteri and S.~Polidoro}, {\em Harnack inequality}, preprint,  (2017).

\bibitem{Fabes1993}
{\sc E.~B. Fabes}, {\em Gaussian upper bounds on fundamental solutions of
  parabolic equations; the method of {N}ash}, in Dirichlet forms ({V}arenna,
  1992), vol.~1563 of Lecture Notes in Math., Springer, Berlin, 1993,
  pp.~1--20.

\bibitem{GarofaloLanconelli}
{\sc N.~Garofalo and E.~Lanconelli}, {\em Level sets of the fundamental
  solution and {H}arnack inequality for degenerate equations of {K}olmogorov
  type}, Trans. Amer. Math. Soc., 321 (1990), pp.~775--792.

\bibitem{4authors}
{\sc F.~Golse, C.~Imbert, C.~Mouhot, and A.~Vasseur}, {\em Harnack inequality
  for kinetic {F}okker-{P}lanck equations with rough coefficients and
  application to the {L}andau equation}, preprint, arXiv.org:1607.08068,
  (2016).

\bibitem{GolseVasseur}
{\sc F.~Golse and A.~Vasseur}, {\em H\"older regularity for hypoelliptic
  kinetic equations with rough diffusion coefficients}, preprint, Version 2,
  (2015).

\bibitem{HobsonRogers}
{\sc D.~G. Hobson and L.~C.~G. Rogers}, {\em Complete models with stochastic
  volatility}, Math. Finance, 8 (1998), pp.~27--48.

\bibitem{Hormander}
{\sc L.~H{\"o}rmander}, {\em Hypoelliptic second order differential equations},
  Acta Math., 119 (1967), pp.~147--171.

\bibitem{Il'in}
{\sc A.~M. Il{\cprime}in}, {\em On a class of ultraparabolic equations}, Dokl.
  Akad. Nauk SSSR, 159 (1964), pp.~1214--1217.

\bibitem{ImbertMouhot}
{\sc C.~Imbert and C.~Mouhot}, {\em H\"older continuity of solutions to
  hypoelliptic equations with bounded measurable coefficients}, preprint,
  arXiv.org:1505.04608, Version 5,  (2015).

\bibitem{Kupcov5}
{\sc L.~P. Kupcov}, {\em Mean value theorem and a maximum principle for
  {K}olmogorov's equation}, Mathematical notes of the Academy of Sciences of
  the USSR, 15 (1974), pp.~280--286.

\bibitem{Kupcov4}
\leavevmode\vrule height 2pt depth -1.6pt width 23pt, {\em On parabolic means},
  Dokl. Akad. Nauk SSSR, 252 (1980), pp.~296--301.

\bibitem{LanconelliPascucci2016}
{\sc A.~Lanconelli and A.~Pascucci}, {\em Nash estimates and upper bounds for
  non-homogeneous {K}olmogorov equations}, to appear in Potential Anal. DOI
  10.1007/s11118-017-9622-1,  (2017).

\bibitem{LanconelliPolidoro}
{\sc E.~Lanconelli and S.~Polidoro}, {\em On a class of hypoelliptic evolution
  operators}, Rend. Sem. Mat. Univ. Politec. Torino, 52 (1994), pp.~29--63.

\bibitem{Langevin}
{\sc P.~Langevin}, {\em On the theory of {B}rownian motion}, C. R. Acad. Sci.
  (Paris), 146 (1908), pp.~530--533.

\bibitem{Lions1}
{\sc P.-L. Lions}, {\em On {B}oltzmann and {L}andau equations}, Philos. Trans.
  Roy. Soc. London Ser. A, 346 (1994), pp.~191--204.

\bibitem{Moser3}
{\sc J.~Moser}, {\em A {H}arnack inequality for parabolic differential
  equations}, Comm. Pure Appl. Math., 17 (1964), pp.~101--134.

\bibitem{Moser3bis}
\leavevmode\vrule height 2pt depth -1.6pt width 23pt, {\em Correction to: ``{A}
  {H}arnack inequality for parabolic differential equations''}, Comm. Pure
  Appl. Math., 20 (1967), pp.~231--236.

\bibitem{Nash}
{\sc J.~Nash}, {\em Continuity of solutions of parabolic and elliptic
  equations}, Amer. J. Math., 80 (1958), pp.~931--954.

\bibitem{pascuccibook}
{\sc A.~Pascucci}, {\em P{DE} and martingale methods in option pricing}, vol.~2
  of Bocconi \& Springer Series, Springer, Milan; Bocconi University Press,
  Milan, 2011.

\bibitem{PP2003}
{\sc A.~Pascucci and S.~Polidoro}, {\em A {G}aussian upper bound for the
  fundamental solutions of a class of ultraparabolic equations}, J. Math. Anal.
  Appl., 282 (2003), pp.~396--409.

\bibitem{PP2004}
\leavevmode\vrule height 2pt depth -1.6pt width 23pt, {\em The {M}oser's
  iterative method for a class of ultraparabolic equations}, Commun. Contemp.
  Math., 6 (2004), pp.~395--417.

\bibitem{Peszek}
{\sc R.~Peszek}, {\em {PDE} models for pricing stocks and options with memory
  feedback}, Appl. Math. Finance, 2 (1995), pp.~211--223.

\bibitem{Polidoro2}
{\sc S.~Polidoro}, {\em On a class of ultraparabolic operators of
  {K}olmogorov-{F}okker-{P}lanck type}, Matematiche (Catania), 49 (1994),
  pp.~53--105.

\bibitem{Polidoro1}
\leavevmode\vrule height 2pt depth -1.6pt width 23pt, {\em A global lower bound
  for the fundamental solution of {K}olmogorov-{F}okker-{P}lanck equations},
  Arch. Rational Mech. Anal., 137 (1997), pp.~321--340.

\bibitem{Risken}
{\sc H.~Risken}, {\em The {F}okker-{P}lanck equation: Methods of solution and
  applications}, Springer-Verlag, Berlin, second~ed., 1989.

\bibitem{Sonin}
{\sc I.~M. Sonin}, {\em A class of degenerate diffusion processes}, Teor.
  Verojatnost. i Primenen, 12 (1967), pp.~540--547.

\bibitem{WangZhang2011}
{\sc W.~Wang and L.~Zhang}, {\em The {$C^\alpha$} regularity of weak solutions
  of ultraparabolic equations}, Discrete Contin. Dyn. Syst., 29 (2011),
  pp.~1261--1275.

\bibitem{Weber}
{\sc M.~Weber}, {\em The fundamental solution of a degenerate partial
  differential equation of parabolic type}, Trans. Amer. Math. Soc., 71 (1951),
  pp.~24--37.

\end{thebibliography}

\end{document}